\documentstyle[amscd,amssymb,verbatim]{amsart}

\theoremstyle{plain}
 
\newtheorem{Lem}{Lemma}
\newtheorem{Prop}{Proposition}
\newtheorem{Def}{Definition}
\newtheorem{Con}{Conjecture}
\newtheorem{Rem}{Remark}

\newcommand{\g}{{\frak g}}
\newcommand{\GG}{\Gamma}
\newcommand{\dd}{\partial}
\newcommand{\LL}{\Lambda^.({\frak g})}
\begin{document}

 \begin{abstract} We introduce a formalism which allows us to formulate a version of Fukaya category in presence of curves of higher genus.
\end{abstract}

\title[Fukaya category  \\
with curves of higher genus]
{Fukaya category \\
with curves of higher genus}

\date{} 
\author{Michael Movshev} 
\address{Department
of Mathematics \\ 
University at Stony Brook \\
Stony Brook, NY} 
\email{mmoushev@@math.sunysb.edu}

\dedicatory{ To the memory of  P.K.}


\maketitle

\section{Stringy category with one object}
\subsection{Frobenius algebras}
Suppose $A$ is a finite dimensional Frobenius algebra, i.e. $\Bbb{Z}$-graded  algebra with a nondegenerate super-symmetric bilinear form $(.,.)$, which satisfies $(ab,c)=(a,bc)$. Assume that the form defines a symmetric  isomorphism $A^*=A[n]$, i.e. we have a nondegenerate (graded) symmetric  pairing $A_i \otimes A_{n-i} \rightarrow k$. We are going to study some generalizations of Frobenius algebras.

 Fix a $\Bbb{Z}$-graded linear space $A$ and a bilinear symmetric  form $(.,.)$ as in the previous paragraph. Denote the grading of $a \in A_i$ by $\tilde a$.  We are willing to treat a pair  $( A[1],(.,.)) $ as a symplectic space. The functor $[n]$, $n \in {\Bbb Z}$ is a shift of the grading: $A[n]^i=A^{i+n}$.  This setup  brings us into a context of formal non-commutative symplectic geometry \cite{K}. In \cite{K} the author introduces three Lie algebras $l_n,a_n,c_n$ for  $2n$ dimensional symplectic  vector space. we shall need an graded version of algebra $a_n$ which we denote by $\tilde{\cal A}$.

In more details : let $T( A[1]^*)$ be a tensor algebra generated by a vector space $A[1]^*$. Denote $D(A^*)$ a super-Lie algebra of super-derivations of $T( A[1]^*)$. Denote the super-commutator by  $[.,.]$. The biliniar form  $(.,.)$ defines an element $\omega$ in $ ( A[1]^*)^{\otimes 2}$, which by construction is in $[T( A[1]^*),T( A[1]^*)]$. Consider a subalgebra $\tilde{\cal A} (A^*)$ of  $D(A^*)$, which preserves $\omega$. The algebra $\cal A$$(A^*)$ carries a grading induced from the grading on $T( A[1]^*)$.

Suppose $A$ is  a Frobenius algebra. Then (assuming that $dim A < \infty$) the structure constants provide us with a tensor
\begin{equation} \label{E:qsd}
 \mu \in   Hom(A[1]^*, A[1]^* \otimes A[1]^*)
\end{equation}
 The associativity condition  translates into $[\mu, \mu ]=0$. Frobenius condition  $(ab,c)$\\$=$$(a,bc)$ translates into  $\mu  \in \tilde{\cal A}^2 \subset Hom(A[1]^*, A[1]^* \otimes A[1]^*)$.

\subsection{$A^{\infty}$ Frobenius algebras}
A well known generalization of a Frobenius algebra is a tensor of a mixed degree   $M=\sum_{i=1}^{\infty} \mu_i$, $deg(\mu_i)=i$ ($i$ is odd) satisfying an equation
\begin{equation} \label{E:poi}
 [M,M]=0
\end{equation}
We shall call  this tensor an   $A^{\infty}$ Frobenius algebra on linear space $A$ . 

The functor of reduction of grading modulo two turns $(T( A[1]^*),M) $ into a differential super algebra.   By definition morphisms between  $(T( A_1[1]^*),M_1) $ and $(T( A_2[1]^*),M_2) $ is a morphism between corresponding differential superalgebras.
\begin{Prop}
Equivalence classes of deformations of solution of equation  \ref{E:poi} are isomorphic to 
cohomology of the complex $\cal A^.$ with a differential being  a commutator with element $M$.  
\end{Prop}

\begin{Rem}
The complex $\tilde{\cal A^.}$ with differential $\mu$ as in \ref{E:qsd} is isomorphic to complex computing cyclic cohomology of algebra $A$.
\end{Rem}
 
\subsection{Alternative definition of $\tilde{\cal A}$}

The following definition is borrowed from \cite{K}.

 Let $T(A[1])$ tensor  algebra over graded linear space $A[1]$. A grading $deg$ on homogenious  generators $a_i \in A[1]$ is equal to $\tilde a_i -1$ Denote ${\cal A}$ the linear space $(T(A[1])/[T(A[1]),T(A[1])])[n-2]$. For every homogeneous element $a \in (A^s[1])^*=A^{-s}[-1] $ one can define two linear maps $\frac{\dd}{\dd\overset{\leftarrow}{ a}}:  {\cal A} \rightarrow T(A[1])[n-2-s]$, $\frac{\dd}{\dd \overset{\to}{a}}:  {\cal A} \rightarrow T(A[1])[n-2-s]$ by the rule:
\begin{equation}\notag
\frac{\dd H}{\dd \overset{\leftarrow}{a}}=\sum_{i=1}^{n}(1\otimes \dots \otimes 1 \otimes a) (\tau)^i \tilde H
\end{equation}
\begin{equation}\notag
\frac{\dd H}{\dd \overset{\rightarrow}{a}}=\sum_{i=1}^{n}(a \otimes 1 \otimes \dots \otimes 1 ) (\tau)^i \tilde H
\end{equation}
In this formula $\tilde H$ is a homogeneous representative of $H$ in  $T(A[1]^*)$. The degree of  $\tilde H$ is $n$. The transformation $\tau$ is a generator of the group ${\Bbb Z}_n \subset S_n$ acting on $n$-th tensor power of graded vector space $A[1]$.

Since we have a symmetric bilinear form $(.,.)$ one can utilize it to combine all $\frac{\dd}{\dd \overset{\to}{ a}}$ into universal map 
\begin{equation}\notag
{\cal A} \rightarrow Hom(A[1],T(A[1]))
\end{equation}
We can identify $A$ and $A^*$ using the bilinear form. Then one can check that the image  of this map is precisely equal to $\tilde {\cal A}$ and the kernel is one-dimensional and generated by the constants. 

{\bf Assumption} From now on we  do not assume  that the bilinear form is invertible.

One can describe the  Lie algebra structure in terms of linear space ${\cal A}$.  Fix a homogeneous basis $a_1 , \dots ,a_n$ the vector space $A[1]$. The commutator of elements $H$ and $G$ is given by the formula:

\begin{equation}
 [H,G]=\sum_{i,j=1}^n \frac{\dd H}{\dd \overset{\leftarrow}{a}_i}(a_i,a_j)  \frac{\dd G}{\dd \overset{\rightarrow}{a}_j}
\end{equation}
According to \cite{K} this is a graded Lie algebra.

\subsection{Co-commutator}
In this section we shall describe a new stricture: a map $\lambda:{\cal A} \rightarrow ({\cal A} \otimes {\cal A})$
Suppose we are given an element $H=a_{i_1}\dots a_{i_k} \in {\cal A} $ which is a cyclic word in $a_1 , \dots ,a_n$  and we would like to compute $\lambda(H)$. 

Denote $H_{l,m}$ the tensor product
\begin{equation}\notag
\pm (a_{i_l},a_{i_m}) \times  a_{i_{m+1}} \dots a_{i_{l-1}} \otimes a_{i_{l+1}} \dots a_{i_{m-1}}
\end{equation}
\begin{equation}
\pm=deg( a_{i_{l+1}} \dots a_{i_{m-1}})deg(a_{i_l})
\end{equation}
of two  cyclic words.

Then
\begin{equation}\notag
\lambda (H)= \sum H_{l,m}
\end{equation}
In the last formula $l < m$ in a natural cyclic order and the summation is taken over all possible $l,m$.
\begin{Prop}
a) $\lambda$ maps ${\cal A}$ into $\Lambda^2({\cal A})$which decreas the grading on $2n-4$.\\
\\
b)The dual map $\lambda^*:\Lambda^2({\cal A}^*) \rightarrow {\cal A}^*$ makes $ {\cal A}^*$  a Lie algebra.\\
\\
c) $\lambda$ is a one-cocycle of ${\cal A}$ in $\Lambda^2({\cal A})$.\\
\\
d) The composition $[.,.] \circ \lambda: {\cal A} \rightarrow {\cal A}$ is equal to zero.
\end{Prop} 

\begin{Rem}
Statements a)-c) of the last proposition can be rephrased as ${\cal A}$ is a Lie bialgebra. Property $d$ doesn't hold for all Lie bialgebras. In general such composition is a derivation and a co-derivation.
\end{Rem}

We would like to explore dependance of ${\cal A}(A,(.,.))$ on graded linear inner product  space $A$. Suppose $(A_1,(.,.))$  is a graded linear inner product subspace of $(A_2, (.,.))$ .
\begin{Prop}
There is a natural inclusion of bialgebras 
\begin{equation}
{\cal A}(A_1,(.,.)_1)\to {\cal A}(A_2,(.,.)_2)
\end{equation}
\end{Prop}
\begin{pf}
The proof is a direct corrolary of definitions.
\end{pf}   

\subsection{Stringy categories with one object}
To clarify the essence of our construction we are going to replace Lie bialgebra ${\cal A}^.$ by an arbitrary graded Lie bialgebra $\g$ with a bracket $[.,.]$ and a  co-commutator $\lambda$.

 The exterior algebra $\Lambda^.(\g)$  is by definition is a summetric algebra $S^.(\g[1])$.

The algebra $\Lambda^.(\g)$ carries tree structures. The first one is a structure of  super-Poisson algebra, coming from interpretation of  $\LL$ as an algebra of left invariant polyvector fields on a supergroup $G$.

\begin{Rem}\label{R:wdg}
The reason of a great care we write all the structures is that one has more than one definition of exterior algebra in the graded category. The altermative definiton of $\Lambda (\g)$ is a quotient of free algebra over $\g$ by an ideal generated by relations $a_1a_2=-a_2a_1, a_1,a_2 \in \g ^{even}$  $b_1b_2=b_2b_1, b_1,b_2 \in \g ^{odd}$ and  $ab=-ba, a \in \g ^{even},b \in \g ^{odd}$. The last equation has a dramatic consequnce: the bracket which one can build doesn't define a graded Lie algebra (but rather its generalization). 
\end{Rem}

 Choose a basis $a_1,\dots,a_n ,\dots,$ in  $\g^{even}$ and $b_1,\dots,b_n ,\dots,$ in  $\g^{odd}$ . An  element $A=a_{r_1}  \dots   a_{r_n}b_{p_1}\dots b_{p_k}$ has a grading
\begin{equation}\label{E:odt}
 Deg(A)=\sum_{t=1}^n (deg(a_{r_s})+1)+\sum_{t=1}^k (deg(b_{p_t})+1)
\end{equation}
 and cohomological grading $\tilde A$ .

    Let $A=a_{r_1}   \dots   a_{r_e}b_{p_1} \dots b_{p_f}$ and  $B=a_{s_1}   \dots   a_{s_g}b_{t_1} \dots b_{t_h}$. The Poisson bracket in  $\LL$ which we denote by by $\{.,.\}$ is given by the formula:
\begin{align}
\{A,B \}=\sum_{i,j} \frac{\dd A}{\dd \overset{\leftarrow}{a}_i}[a_i,a_j]  \frac{\dd B }{\dd \overset{\rightarrow}{a}_j}+\sum_{i,j} \frac{\dd A}{\dd b_i}[b_i,b_j]  \frac{\dd B }{\dd b_j}+ \notag \\
\sum_{i,j} \frac{\dd A}{\dd \overset{\leftarrow}{a}_i}[a_i,b_j]  \frac{\dd B }{\dd b_j}+\sum_{i,j} \frac{\dd A}{\dd b_i}[b_i,a_j]  \frac{\dd B }{\dd \overset{\rightarrow}{a}_j}
\end{align}

  The second structure is a homological differential $\dd$.

\begin{align}
& \dd(A)= - \sum_{i<j }[a_i,a_j] \frac{\dd^2 A}{\dd \overset{\to}{a}_i \dd \overset{\to}{a}_j} + \sum_{i<j }[b_i,b_j] \frac{\dd^2 A}{\dd b_i \dd b_j} + \sum_{i,j }[a_i,b_j] \frac{\dd^2 A}{\dd \overset{\to}{a}_i \dd b_j}
\end{align}

 These structures interact  in a natural way.

\begin{Prop}
The maps    $\{.,.\}: \LL \otimes \LL \rightarrow \LL$  and $ \dd: \LL \rightarrow \LL$ have cohomological equal to $-1$. Also  $Deg(\{.,.\})=Deg( \dd)=-1$ and these operations satisfy  an equation: 
\begin{equation} \notag
(-1)^{Deg(A) }\{A,B\}=(\dd A)B - \dd(AB)+(-1)^{Deg(A)}A(\dd B)
\end{equation}

From this equation it follows that $\dd$ is a differentiation of the bracket $\{.,.\}$.

The bracket $\{.,.\}$ satisfies a Lie algebra relations:
\begin{align}
&\{A,B\}=(-1)^{Deg(A)-1)(Deg(B)-1)+1} \{ B,A \} \notag \\
&\{A,\{ B,C \} \}= \{ \{A,B \},C \}+(-1)^{Deg(A)-1)(Deg(B)-1)} \{B,\{A,C \} \} \notag
\end{align}
A precise formulation that we have a Poisson algebra is as follows:
\begin{equation}
\{A, BC \}=  \{A,B\}C+(-1)^{(Deg(A)-1)Deg(B)} B \{A,C \}
\end{equation}
\end{Prop}

\begin{pf}We shall work out only purely even case of the first statement. The general proof is left as exercise to the reader. The rest of the statements are fairly standard and the proof will be ommited. A refference for a sign rules which may be used in the proof of the formulas is \cite{F}
 
By definition 
\begin{align}
&\dd ((a_1   \dots   a_i)   (a_{i+1}   \dots   a_k))\overset{def}{=}\notag \\
&\sum_{s<t}(-1)^{s+t-1}[a_s,a_t] (a_1   \dots \Hat a_s \dots \Hat a_t \dots   a_i)   (a_{i+1}   \dots   a_k)+\notag\\
&\sum(-1)^{s+t-1}[a_s,a_t] (a_1   \dots \Hat a_s \dots    a_i)   (a_{i+1}  \dots  \Hat a_t \dots   a_k)+\notag \\
&\sum_{s<t}(-1)^{s+t-1}[a_s,a_t] (a_1    \dots    a_i)   (a_{i+1} \dots \Hat a_s  \dots  \Hat a_t \dots   a_k)= \notag \\
&=\dd( a_1   \dots   a_i )   a_{i+1}   \dots   a_k+\notag \\
&+(-1)^{deg(a_1   \dots   a_i)+1}\{a_1   \dots   a_i,a_{i+1}   \dots   a_k\}+\notag \\
&+(-1)^{deg(a_1   \dots   a_i)}a_1   \dots   a_i   \dd( a_{i+1}   \dots   a_k) \notag
\end{align}
Proof that $\dd$ is a differentiation:
\begin{align}
& \dd \{ a,b \}=(-1)^a(\dd (( \dd a)b)+(-1)^a \dd (a \dd b))= \notag \\
&= \{ \dd a,b \}- \dd a \dd b+ \dd a \dd b +(-1)^{a+1}\{a,\dd b \})= \notag\\
&=\{ \dd a,b \}+(-1)^{a+1}\{a,\dd b \} \notag
\end{align}
\end{pf}
The last structure is a cohomological differential $d$. It is a differentiation of the exterior algebra. Let 
\begin{align}
&\lambda_{0,0}:\g^{even} \rightarrow \Lambda^2(\g^{even})\notag \\
&\lambda_{1,1}:\g^{even} \rightarrow S^2(\g^{odd})\notag \\
&\lambda_{0,1}:\g^{odd} \rightarrow \g^{even}\otimes \g^{odd}\notag \\
&\lambda_{1,0}:\g^{odd} \rightarrow \g^{odd} \otimes  \g^{even}\notag 
\end{align}
be components of the co-commutator
 On a generator $a \in \g^{even}$ it is defined as $d(a)=\lambda_{0,0}(a)-\lambda_{1,1}(a)$,on generator $b \in \g^{odd}$ the differential is $d(b)=\lambda_{0,1}(b)-\lambda_{1,0}(b)$
\begin{Rem}
One can try to take $d(x)=\lambda(x)$ but in this case the image will land in a wrong algebra (the one defined in remark \ref{R:wdg})
It is easy to see that $d^2=0$ and it can be interpreted as a cohomological differential of the dual Lie algebra.
\end{Rem}

\begin{Prop}
$d$ is a differential on the Poisson algebra $\LL$ of  cohomological degree $1$. In case of algebra ${\cal A}$ its degree $Deg(d)$ is $4-2n$.
\end{Prop}
\begin{pf}
We prove the statement in the simplest case when algebra $\g$ is concentrated in degree zero.
By constriction $d$ is a differential of the exterior algebra.
One needs to check that $d$ is a differentiation of the bracket $\{A,B\}$.It can be proved by induction on $Deg(A)+Deg(B)-2$.
Indeed when $Deg(A)=Deg(B)=1$ then the condition  $ d \{ A,B \}= \{ d A,B \}+(-1)^{A+1}\{A,d B \}$ is a restatement that $\lambda$ is a one-cocycle.
Suppose that $Deg(A)+Deg(B)+Deg(C)=n$ and $Deg(B),Deg(C) \geq 1$. In the third identity we use inductive assumptions:
\begin{align}
&d \{ BC,A \}=d (B\{ C,A \}+(-1)^{(Deg(A)-1)Deg(C)}\{ B,A \}C)= \notag \\
&=d(B)\{ C,A \}+(-1)^{Deg(B)} B( \{d( C),A \}+ (-1)^{Deg(C)+1}\{ C,d(A) \})+\notag \\
& +(-1)^{(Deg(A)-1)Deg(C)}( \{ d(B),A \}+(-1)^{Deg(B)+1} \{ B,d(A) \})C+\notag \\
&+(-1)^{(Deg(A)-1)Deg(C)+Deg(B)+Deg(A)-1} \{ B,A \}d(C) =\notag \\
&=\{ d(B)C+(-1)^{Deg(B)} B d(C),A \}+(-1)^{Deg(B)+Deg(C)+1}\{ BC,d(A) \} \notag
\end{align}
\end{pf} 

\begin{Prop}
$d \circ \dd +\dd\circ d=\psi$ The map $\psi$ is a derivation of Poisson algebra $\LL$. On generating linear space $\g$ it is equal to composition $[.,.] \circ \lambda$.
\end{Prop}
\begin{pf}
To shorten the exposition we prove only for algebra $\g$ which has only zero graded component.
\begin{align}
&\dd\circ d(AB)+d\circ \dd(AB)=\notag \\
&=\dd(d(A)B+(-1)^{Deg(A)}Ad(B))+d(\dd(A)B+(-1)^{Deg(A)+1}\{A,B\}+\notag \\
&+(-1)^{Deg(A)}A\dd(B))=\notag \\
&=\dd\circ d(A)B+(-1)^{Deg(A)} \{ d(A),B\}+(-1)^{Deg(A)+1}d(A)\dd (B)+ \notag 
\end{align}
\begin{align}
& +(-1)^{Deg(A)} \dd (A)d(B)+(-1)^{2Deg(A)-1}\{A,d(B)\}+A \dd\circ d(B)+d\circ \dd (A)B+\notag \\ 
&+(-1)^{Deg(A)-1} \dd (A)d(B) +(-1)^{Deg(A)+1}d \{ A,B \}+(-1)^{Deg(A)}d(A)\dd(B)+\notag \\
&+(-1)^{2Deg(A)}A d \circ \dd(B)+(-1)^{Deg(A)+1}d\{A,B\}+ \notag \\
&+(-1)^{Deg(A)}\{d(A),B\}-\{A,d(B)\}=\notag\\
&=(\dd\circ d(A)+d\circ \dd(A))B+A(\dd\circ d(B)+d\circ \dd(B)) \notag
\end{align}
The second statement of the proposition is obvious.
\end{pf}

The zero modes of operator $\psi$ form a Poisson subalgebra $\LL_0$ of $\LL$. It is equipped with the action of two anticommuting operators $\dd$ and $d$.

The element $M \in \g_0=\g \cap \LL_0$ can be thought of as an element of $\LL$. Suppose $deg(M)+1=Deg(M)$ is even
 \begin{equation}\label{E:tyu}
[M,M]=0
\end{equation}
This equation materializes in Fukaya category as as equations of $A^{\infty}$ category.

A natural generalization of equations \ref{E:tyu} is
\begin{equation} \label{E:oiy}
d R=\frac{1}{2}\{R,R \}
\end{equation}
Here we replace $M \in \g_0 \subset \LL_0$ by a generic element $R=\sum_{i=1}^{\infty} R_i$ ($R_i \in \Lambda^i({\frak g})$ ) , with all degree $Deg(\tilde R_i)$ are being even (don't forget about a shift of the grading \ref{E:odt}) . The element $M =R_1$ is a solution of \ref{E:tyu}. We are going to call equation \ref{E:oiy} master equation at a tree level.It should appear in a version of Fukaya category where holomorphic disks are replased by general genus zero curves with boundary.

We shall need to make a further generalization of equation \ref{E:oiy}. Let $\alpha$ be a formal variable. Introduce an operator

\begin{equation}\label{E:cft}
\delta_{\alpha}=d +\alpha \dd
\end{equation}

Our  ultimate goal is the equation 
\begin{equation} \label{E:svy}
\delta_{\alpha}R(\alpha)=\frac{1}{2}\{R(\alpha),R(\alpha) \}
\end{equation} 
Here an even  element $R(\alpha)$ is a formal power series $R(\alpha)=\sum_{g=0} R_g\alpha^g \in \LL_0 [\alpha]$.
This element  should correspond to structure constants of a a version of Fukaya category with no restriction on genus of curves used in its  definition.

\begin{Rem}
Observe that setting $\alpha$ to zero we degenerate equation \ref{E:oiy} to \ref {E:cft}.
\end{Rem}

\begin{Lem}
For any (non homogenious) element $A=\sum A_{2i}$, such that $Deg(A_{2i})=2i$, the following identity holds:
\begin{equation}
\delta_{\alpha}(A^n)=nA^{n-1}\delta_{\alpha}(A)-\frac{n(n-1)}{2}\alpha A^{n-2}\{A,A\}
\end{equation} 
\end{Lem}\label{L:ewq}
\begin{Prop}
Element $exp(\frac{1}{\alpha}R(\alpha))$ satisfies an equation
\begin{equation}
\delta_{\alpha}exp(\frac{1}{\alpha}R(\alpha))=0
\end{equation}
\end{Prop}
\begin{pf}
\begin{align}
&\delta_{\alpha}\sum_{n=0}^{\infty}\frac{R(\alpha)^n}{\alpha^n n!} = \notag \\
&= \sum_{n=0}^{\infty}\frac{nR(\alpha)^{n-1}}{\alpha^n n!}\delta_{\alpha}R(\alpha)- \sum_{n=0}^{\infty}\frac{n(n-1)\alpha R(\alpha)^{n-2}}{2 \alpha^n n!}\{R(\alpha),R(\alpha)\}= \notag \\
&=\frac{1}{\alpha}exp(\frac{1}{\alpha}R(\alpha))(\delta_{\alpha}R(\alpha)-\frac{1}{2}\{R(\alpha),R(\alpha)\})=0
\end{align}
We used in the proof the result of lemma \ref{L:ewq}
\end{pf}

\subsection{BV algebra}

\begin{Def}\label{D:iuo}

By definition  BV-algebra $B$  is a ${\Bbb Z}_2$ graded unital  Poisson algebra with an odd bracket and an odd operator  $\dd$ such that 
\begin{equation}\notag
(-1)^{\tilde a}\{a,b\}=-\dd(ab)+(\dd a)b+(-1)^{\tilde a}a(\dd b) 
\end{equation} 
and $\dd^2=0$
\end{Def}

 The group of invertible elements of algebra $B$ acts on odd operators $D$ whose square is zero by conjugation.

Let us see what is an orbit of $d$. Indeed
\begin{equation}
d_a(x)=a^{-1}\dd(ax)=\dd(x)-a^{-1}\{a,x\}+a^{-1}\dd(a)x
\end{equation}
It is clear that $\dd(1)=0$. In order for $\dd_a(x)$ to be a part of a structure defining a BV-algebra $\dd_a(1)$ must be zero, so $\dd(a)=0$.

It is easy to prove the converse
\begin{Prop}
The  ${\Bbb Z}_2$ graded unital  Poisson algebra from definition \ref{D:iuo}  with a new  odd differential $\dd_a(x)$ is a BV-algebra iff $\dd(a)=0$.  
\end{Prop}
We are looking for a natural equivalence relation in the space of invertible $a$  such that $\dd(a)=0$.

The Lie algebra $\g(B)$ assosiated with BV-algebra $B$ acts on $B$, however the action doesn't commute with $\dd$. Our goal is to modify the action of  $\g(B)$, so that the new action would commute with $\dd$.

\begin{Lem}
The formula
\begin{equation}
A(r)a=\{r,a\}+(-1)^rdra
\end{equation}
defines an action of the Lie algebra $\g(B)$. It commutes with differential: $\dd A(r)=(-1)^{r}A(r)\dd$
\end{Lem}
It is clear that the action is twisted by one-cocycle $\dd$ with walues in $B$.

It is possible to describe the corresponding group action on $B^*$. It is clear that it will be twisted on integral of cocycle $\dd$.
\begin{Prop}
Assume $g=exp(r)$
The following formula defines a cocycle of group $G(B)$ with values in $B$.
\begin{equation}
c(g)=\sum_{n=0}^{\infty} \frac{ad(r)^n(\dd r)}{(n+1)!}
\end{equation}
The multiplicative cocycle will be of the form $c^*(g)=exp(c(g))$
\end{Prop}
\begin{pf}
The proof follows from the identity 
\begin{equation}
\frac{d}{dt}c(exp(rt))=exp(rt)\frac{d}{d\epsilon}c(exp(\epsilon r))\Big|_{\epsilon=0}
\end{equation}
 valid for all additive  one-cocycles and identification $c(exp(rt))\Big|_{t=0}=\dd(r)$
\end{pf}

\begin{Lem}\label {L:qaw}
The following identity holds:
\begin{equation}
\sum _{s=0}^{l}\frac{(s+t)!}{s!}=\frac{(l+t+1)!}{l!(t+1)}
\end{equation}
\end{Lem}

\begin{Lem}\label{L:rfb}
 In any BV-algebra a following identity holds:
\begin{equation}
\dd exp(r)a=exp(r)\dd a+\{c(exp(r)),exp(r)a\}
\end{equation}
\end{Lem}
\begin{pf}
\begin{align}
&\dd exp(r)a=\sum_{k=0}^{\infty}\dd \frac{ad^k(r)a}{k!}= \notag \\
&\sum_{k=0}^{\infty}(\sum_{i=0}^{k-1}\frac{1}{k!}ad^i(r)\{\dd r,ad^{k-1-i}(r)a\}+\frac{1}{k!}ad^k(r)\dd a) \notag \\
&\sum_{k=1}^{\infty}\sum_{i=0}^{k-1}\sum_{j=0}^{i}\frac{i!}{k!j!(i-j)!}\{ad^j(r)\dd r,ad^{k-1-j}(r)a\} +exp(r)\dd a = \notag \\
&\sum_{k=1}^{\infty}\sum_{j=0}^{k-1}\sum_{i=j}^{k-1}\frac{i!}{k!j!(i-j)!}\{ad^j(r)\dd r,ad^{k-1-j}(r)a\} +exp(r)\dd a = \label{E:ois} \\
&\sum_{k=1}^{\infty}\sum_{j=0}^{k-1}\frac{1}{(k-j-1)!(j+1)!}\{ad^j(r)\dd r,ad^{k-1-j}(r)a\} +exp(r)\dd a = \notag \\
&=\{c(exp(r)),exp(r)a\}+exp(r)\dd a \notag
\end{align}
To simplify equation \ref{E:ois} we used results of lemma \ref{L:qaw} 
\end{pf}

\begin{Lem} \label{L:dft}
Cocycle  $c(exp(r)$ satisfies an eqiation $$\dd(c(exp(r))=\frac{1}{2}\{c(exp(r),c(exp(r)\}$$.
\end{Lem}
\begin{pf}
\begin{align}
&\dd c(exp(r))=\sum_{k=1}^{\infty}\sum_{i=0}^{k-1}\sum_{j=0}^{i}\frac{i!}{(k+1)!j!(i-j)!}\{ad^j(r)\dd r,ad^{k-1-j}\dd r\}  = \notag \\
&\sum_{k=1}^{\infty}\sum_{j=0}^{k-1}\frac{(k-j)}{(k+1)}\{\frac{ad^j(r)\dd r}{(j+1)!},\frac{ad^{k-1-j}(r)\dd r}{(k-j)!}\}  = \notag \\
&\frac{1}{2}\sum_{k=1}^{\infty}\sum_{j=0}^{k-1}\{\frac{ad^j(r)\dd r}{(j+1)!},\frac{ad^{k-1-j}(r)\dd r}{(k-j)!}\}  =\frac{1}{2}\{c(exp(r),c(exp(r)\} \notag 
\end{align}
\end{pf}
\begin{Prop}\label{E:htr}
Define an action of group $G(B)$ on $B$ via the formula
\begin{equation}
T(exp(r))a=c^*(exp(r))exp(r)a
\end{equation}
The differential $\dd$ commutes with all operator $T(exp(r))$.
\end{Prop}

\begin{pf}

From lemma \ref{L:rfb} we have:
\begin{equation}
exp(r)\dd exp(-r) a=\dd a-\{c(exp(r)),a\}
\end{equation}
 Lemma \ref{L:dft} implies:
\begin{equation}
exp(-c(exp(r))\dd exp(c(exp(r)) a=\dd a-\{c(exp(r)),a\}
\end{equation}
Proposition follows from it.
\end{pf}

\begin{Def}
We say that two solutions $a_1,a_2 \in B^*$ of master equation are equivalent if there is an element $g \in G(B)$ such that $T(g)a_1=a_2$ 
\end{Def}
\begin{Rem}
It is possible to define a superspace ${\cal B}*$ of invertible elements in $B$. It is also possible to define a supergroup ${\cal G}(B)$ and an action on${\cal B}*$. In addition to it operator $\dd$ can be extended as a vector field on ${\cal B}*$. One can check that  ${\cal G}(B)$ commutes with vector field $\dd$.  
\end{Rem}

\begin{Prop}
a)Infinitesimal classes of deformations $r$ of solution $a$ of master equation are in one-to-one correspondence with even  cohomology classes of the complex $(B,\dd-a^{-1}\{a,.\})$.\\
b) The commutator on cohomology is equal to zero.\\
c) Cohomology groups for different $a$ are (non-canonically) isomorphic. \\
\end{Prop}
\begin{Rem}
The full cohomology groups form a  tangent space to super coset space of solutions of master equation.  
  \end{Rem}
\begin{Def}
We say that a stringy category with one object is given if 

a) we fix a graded vector space $A$ with a symmetric bilinear form $(.,.)$ of degree $n$.

b) in the linear space $\Lambda^.({\cal A}(A[1])$ we fix and element  $R(\alpha)$ of total odd degree .

c) This element satisfies equation \ref{E:svy} 
\end{Def}

{\bf Example.} An ordinary linear catigory with one object is an algebra suppose it is a Frobenuis algebra. Let us see how it becomes a particular case of a stringy category.

 Suppose we have a graded finite dimensional Frobenuis algebra $A$ with structure constants $c^i_{j,k}$. One can rase indeces with the bilinear form and get an element $c^{i,j,k} \in {\cal A}^3(A[1])$.It is easy to see that 
$\dd (c)=d(c)=\{c^{i,j,k}.c^{i,j,k}\}=0$. So we have a stringy category. 

\section{General stringy category}
The relation of stringy category to ordinary $A_{\infty}$ category is the same as a relation of stringy algebras to $A_{\infty}$ algebras.

{\bf Objects of $C$} Assume that objects form a set .

{\bf Morphisms of $C$} For any pair of oblects $L_1,L_2$ , $Hom(L_1,L_2)$ is a graded  finite dimensional linear space. There is a pairing  
\begin{equation} \label {E:dfv}
\omega: Hom(L_1,L_2)^i\otimes  Hom(L_2,L_1)^{n-i} \rightarrow k
\end{equation}

{\bf Compositions of morphisms} This the most interesting part of the definition. For any $k$-tuple of sequences of objects 
\begin{align}
&(L_1(1), \dots, L_{p_1}(1) ),\notag \\
& \dots,\notag\\
& (L_1(k), \dots, L_{p_k}(k))\notag
\end{align}
 one has a one parametric family of  tensors
\begin{align}
&R(\alpha)(L_1(1), \dots, L_{p_1}(1)), \dots , (L_1(k), \dots, L_{p_k}(k)) \in\notag  \\
& (Hom(L_1(1),L_2(1)) \otimes \dots \otimes Hom(L_{p_1}(1), L_{1}(1)) ) \otimes \notag \\
&\dots \notag\\
& \otimes  (Hom(L_1(k),L_2(k)) \otimes \dots \otimes Hom( L_{p_k}(k), L_{1}(k))) \notag \\
\end{align}
All these  $R$-s satisfy certain equations. The equations are generalizations of $A_{\infty}$ equations in  $A_{\infty}$ category.

First, the tensors $R$ are defined up to graded cyclic permutation of \\ $ (L_1(i), \dots, L_{p_i}(i))$ for every $i$. 

Second, the tensors $R$ are graded skew symmetric with respect to permutation of  $ (L_1(i), \dots, L_{p_i}(i))$ and  $ (L_1(j), \dots, L_{p_j}(j))$.

A vector space $\bigoplus_{L_1,L_2 \in Ob(C)} Hom(L_1,L_2)$ is graded and  carries a graded symmetric form which is a direct sum of \ref{E:dfv}.

By the above definitions   the element ${\Bbb R}(\alpha)=\bigoplus R(\alpha)$ (the direct sum of all $R(\alpha)$'s introduced above) belongs to $\Lambda^.({\cal A}(\bigoplus_{L_1,L_2 \in Ob(C)} Hom(L_1,L_2)[1]))$. The set of equations satisfied by ${\Bbb R}$ can be compressed into a  formula:
\begin{equation}\notag
\delta_{\alpha}{\Bbb R}(\alpha)=\frac{1}{2} \{{\Bbb R}(\alpha),{\Bbb R }(\alpha) \}
\end{equation} 
The operator $\delta{\alpha}$ was introduced in equation \ref{E:cft}.

One can get an  interesting refinements of this definition  by incorporating an  additive structure in the category $C$. Suppose that there is a functor 
\begin{equation}\notag
\oplus:Ob(C)\times Ob(C) \rightarrow Ob(C)
\end{equation}
This map is  associative.
On the morphisms we have
\begin{align}
& Hom(K \oplus L,M) = Hom(K,M) \oplus Hom(L,M)\notag \\
& Hom(K,L \oplus M) = Hom(K,L) \oplus Hom(K,M)\notag
\end{align}
The important is that for compositions we have 
\begin{align}
&R(\alpha)(L_1(1), \dots, L_{p_1}(1)), \dots,(L_1(i), \dots, L_j(i) \oplus \tilde L_j(i), \dots, L_{p_i}(i)) \dots ,\notag \\
& (L_1(k), \dots, L_{p_k}(k))= \notag \\
&R(\alpha)(L_1(1), \dots, L_{p_1}(1)), \dots,(L_1(i), \dots, L_j(i), \dots, L_{p_i}(i)) \dots , (L_1(k), \dots, L_{p_k}(k))\notag\\ 
& \oplus \notag \\
&R(\alpha)(L_1(1), \dots, L_{p_1}(1)), \dots,(L_1(i), \dots,\tilde L_j(i), \dots, L_{p_i}(i)) \dots , (L_1(k), \dots, L_{p_k}(k)) \notag
\end{align}

\begin{Prop}
 The functor $\oplus$ induce a structure of algebra on homology of the complex
$(\Lambda({\cal A}(\bigoplus Hom(L_1,L_2)[1])), {\Bbb R})$
\end{Prop}

\begin{Rem}
According to the last proposition we have a structure of an algebra on the tangent space of the moduli of additive stringy categories.
\end{Rem}
\begin{Rem}
The moduli of all additive stringy categories of given dimension  is connected.(Is it a math formulation of smooth topology change?) 
\end{Rem}

\section{Modification of Fukaya category}\label{S:bfp}

We borrow some pieces of the following  definition from \cite{K1}.

Let $M$ be a closed symplectic manifold with $c_1(M)=0$

Denote ${\cal L}(M)$ a space of pairs $x,L$, where $x$ is a point in $M$ and $L$ is a Lagrangian subspace in $T_x(M)$. The space  ${\cal L}(M)$  is fibered over $M$ with fibers Lagrangian Grassmanians. The fundamental group of the fiber is isomorphic to  $\Bbb Z$. The condition $c_1(M)=0$ insures that there is a covering $\tilde {\cal L}(M)$ of  ${\cal L}(M)$ inducing a universal covering on each fiber, commuting with a projection on $M$.

An object of  Fukaya category is a  Lagrangian submanifold with a lifting of the canonical map $L \rightarrow {\cal L}(M)$ to a map $L \rightarrow \tilde {\cal L}M$. Also the Lagrangian submanifold $L$ is equipped with a local system $E$. Denote the object $(L,E)$ by $\cal F$

For subvarieties $L_1,L_2$ intersecting each other transversally at a point $a$ one can define a Maslov index $\mu_a(L_1,L_2) \in {\Bbb Z}$. Notice that 
\begin{equation}\notag
\mu_a(L_1,L_2)+\mu_a(L_2,L_1)=n =\frac{1}{2}dim(M)
\end{equation}

Denote  $E_a$ the fiber of the local system $E$ over the point $a \in L$

By definition 

\begin{equation}\notag
Hom({\cal F}_1,{\cal F}_2)= \bigoplus_{a \in L_1 \cap L_2} Hom((E_1)_a,(E_2)_a) \otimes [a]
\end{equation}
with ${\Bbb Z}$-grading coming from Maslov index.

The category $F(M)$ has an action of the shift functor. Denote its action on the object $\cal F$ by ${\cal F}[1]$. The functor $[1]$ shifts the map $L \rightarrow \tilde {\cal L}(M)$ by a positive generator of the deck transformations.

{\bf Convention}In the following discussion we shall consider a set $X(n)\{1, \dots, n \}$ as a homogeneous space over a group ${\Bbb Z}_n$. . It means that if $i<n$ then  $i+1$ is defined as for ordinary integers. However $n+1=1$.

 Fix $k$ tuples of objects $\{ {\cal F}_i(j)|i \in X(p_j) \}$ with underlying Lagrangian submanifolds meeting each other transversally. 
 The tensor \\
\\
\\
\begin{align}
& \ \ \ \ \ \ \ \ \ \ \ \ \ \ \ \ \ \ \ \ \ \ \ \ \ \ \ \ \ \ \ \  R^k(\alpha)\notag \\
& \ \ \ \ \ \ \ \ \ \ \ \ \ \ \ \ \ \ \ \ \ \ \ \ \ \ \ \ \ \ \ \ \ \   \in \notag \\
&\begin{pmatrix}
&Hom({\cal F}_1(1),{\cal F}_2(1))\otimes  & \dots &\otimes Hom({\cal F}_{p_1}(1),{\cal F}_{1}(1))  \\
& \otimes & \dots & \otimes \\
&Hom({\cal F}_1(k),{\cal F}_2(k))\otimes & \dots &\otimes Hom({\cal F}_{p_k}(k),{\cal F}_{1}(k))\\
\end{pmatrix}\notag \\
\end{align}\\
\\
\\

 will be  given   as a formal power series in variable $\alpha$
\begin{equation}\notag
R^k(\alpha)=\sum_{g=0}^{\infty} R^k_g \alpha^{g}
\end{equation}

 Fix elements $u_i(j) \in Hom({\cal F}_i(j)),{\cal F}_{i+1}(j))$($i \in X(p_j)$)  , which are represented by  elements $t_i(j)) \otimes a_{i}(j) $ such that
\begin{equation}\notag
(a_{i}(j),t_i(j)) \in  L_i(j) \cap L_{i+1}(j), Hom((E_i(j))_{a_i(j)},(E_{i+1}(j))_{a_i(j)} 
\end{equation}

 Fix a curve $X$ of genus $g$ with $k$ boundary components. Decompose  the boundary $\dd X = \bigcup S^1(j)$ into a union of intervals
\begin{equation}\notag
S^1(j)= [\alpha_1(j),\alpha_2(j) ] \cup [ \alpha_2(j),\alpha_3(j) ] \cup \dots  [ \alpha_{p_j}(j),\alpha_1(j) ]
\end{equation} 

By definition $R^k_g$ is equal to
\\
\\
\\
$\begin{aligned}
& R^k_g= \sum
\begin{Sb}
a_{i}(j) \in \\
\qquad \qquad \qquad \in L_i(j) \cap  L_{i+1}(j),\\
\end{Sb}
 \begin{matrix}
      & a_1(1)\otimes    & \dots & \otimes a_{p_1}(1)&  \\
	& \otimes  & \dots &  \otimes  \\
       &  a_1(k)\otimes & \dots &\otimes a_{p_k}(k) \\
  \end{matrix} \times \\
& \times C
 \begin{pmatrix}
      & a_1(1),   & \dots & a_{p_1}(1),&  \\
	&   & \dots &    \\
       &  a_1(k), & \dots & a_{p_k}(k) \\
  \end{pmatrix}
\end{aligned}$\\
\\
Where\\
\\

$\begin{aligned}
& C
 \begin{pmatrix}
      & a_1(1),   & \dots & a_{p_1}(1),&  \\
	&   & \dots &    \\
       &  a_1(k), & \dots & a_{p_k}(k) \\
  \end{pmatrix}
= \\
& \ \ \ \ \ \ \ \ \ \ \ \ \ \ \ \ \ \ \ \ \ \ \ \ \ \ \ \ = \sum_{\phi}
 \pm
 e^{2\pi i \int \phi^* \omega}
\Sigma(1) \otimes \dots \otimes \Sigma(k)
\end{aligned}$ \\
The tensor $\Sigma(j)$ will be defined a little later. In the formula we sum over holomorphic maps $\phi: X \rightarrow
M$ up to holomorphic  equivalence,
with the following conditions along the boundary:
there are $p_{j}$ points $\alpha_i(j) \in S^1(j) \subset \dd X$
such that $\phi(\alpha_i(j)) = a_i(j)$ and
$\phi([a_i(j),a_{i+1}(j)]) \subset L_i(j)$.

The tensor $\Sigma(j)$ is in $Hom(E_1(j),E_2(j))_{a_1(j)}\otimes \dots \otimes Hom(E_{p_j}(j),E_1(j))_{a_{p_j}(j)}$ and is a tensor product of monodromies
\begin{equation}\notag
\Sigma_i(j):E_i(j)_{a_i(j)} \rightarrow E_{i+1}(j)_{a_{i+1}(j)}
\end{equation}
along the path $\phi([a_i(j),a_{i+1}(j)])\in L_i(j)$.
\begin{Con}
Elements $R^k(\alpha)$ form structure constants of a stringy category.
\end{Con}
\begin{Con}
Diffrent choices of almost complex strucure $J_1$ and $J_2$ correspond to different elements ${\Bbb R}_1$  ${\Bbb R}_2$. We conjecture that $exp({\Bbb R}_1)$ and $exp({\Bbb R}_2)$ are $G(\Lambda^*({\cal A}))$ equivalent  with respect to the action described in equation \ref{E:htr}. 
\end{Con}

The idea of the proof shouldn't be much different from original Fukaya proof \cite{Fuk}. Roughly speaking the equations of stringy category correspond to the boundary of compactified  moduli of certain  stable open holomorphic  curves in symplectic manifold $M$.

To prove this conjecture one has to have a detailed understanding of certain moduli spaces of open holomorphic curves.

\subsection{Moduli of holomorphic curves with  boundary}
Let us fix the  notations.
Fix $k$ tuples of Lagrangian submanifolds $\{  L_i(j)|i \in X(p_j) \}$($1 \leq j \leq k$)  meeting each other transversally. Sets $X(p_j)$ were introduced in   \ref{S:bfp}

Let $M$ be a manifold with a symplectic form $\omega$ and almost complex structure $J$ compatible with $\omega$ . Let $L_1(j), \dots, L_{p_j}(j)$ be a finite family of  Lagrangian submanifolds. Let $X$ be a holomorphic curve of genus $g$ with $k$ boundary components $S_1,\dots,S_k$. we fix on each boundary component $S_i$ distinct points $x_1(j), \dots ,x_{p_j}(j)$. Fix a homotopy class of a map $[ f ]$ of $X$ into $M$ such that the boundary of $X$ maps into the inion of$L_1(j), \dots, L_{p_j}(j)$ defined up to the action mapping class group fixing marked point on the boundary of $X$ . Moreover the marked points on the boundary map onto intersection points of the Lagrangian submanifolds.  

Denote the moduli of holomorphic map in this homotopy class by $M(f)$.

\begin{Con}

The moduli space  $M(f)$  admits a  compactification to orientable manifold  $\tilde M(f)$   with corners. Moreover one can choose a  orientation uniformly for all  $\tilde M(f)$.
\end{Con}

 We are  interested in the structure of codimension  one boundary components of $\tilde M(f)$ .

Let $t$ be a simple path in $X$ with no self-intersections. The end points of the path belong to the boundary of $X$ and lie off the marked points. 
\begin{Con}
There by the above
e is one-to-one correspondence between isotopy classes of such paths defined up to the action mapping class group fixing marked point on the boundary of $X$  an codimension one components of the boundary of $\tilde M(f)$ .
\end{Con}

 \section {Cohomology of the complex $(\Lambda^.({\cal A}(A[1]),\delta_{\alpha})$}

The cohomology of the complex 
\begin{equation}\label{E:oqx}
(\Lambda^.({\cal A}(A[1]),\delta_{\alpha})
\end{equation}
 when $\alpha = \infty$ has been studied in \cite{K}. There the author uses certain CW-comlexes whose cohomology turns to be closely realted to  $H^.(\Lambda^.({\cal A}(A[1]),\delta_{\alpha})$.

We do not know modify  constructions from  \cite{K} to make them fit in our setting. Instead we introduce followimg \cite{M} a certain small category whose nerve will play the role of mentioned CW-complexes.

 But first we shall elaborate on the structure of the complex $(\Lambda^.({\cal A}(A[1]),\delta_{\alpha})$. Intoduce a complex $(C^.,\delta_{\alpha})$ quasiisimorphic to  $(\Lambda^.({\cal A}(A[1]),\delta_{\alpha})$. The group $C^k$ is  a space of invariants of ${\frak sp}(\infty)\oplus{\frak so}(\infty)$ action on  $(\Lambda^k({\cal A}(A[1]),\delta_{\alpha})$. The groups of the  complex were described in \cite{K}. In the degree $k$ it is generated by linear combinations of fat graphs with $k$ verices (graphs might be disconnected). Introduce a finer grading on this complex. The group $C^{k,j}$ is generated by graphs with $j$ edges and $k$ vertices. The action of the differential is shown on the next picture 

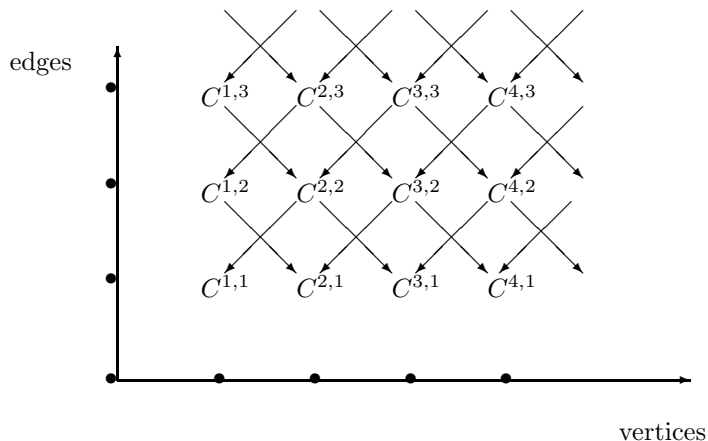
\setlength{\unitlength}{3947sp}
\begin{picture}(4287,2712)(526,-1936)
\thinlines
\put(1801,-1590){\makebox(0,0)[lb]{$\bullet$}}
\put(2401,-1590){\makebox(0,0)[lb]{$\bullet$}}
\put(3001,-1590){\makebox(0,0)[lb]{$\bullet$}}
\put(3601,-1590){\makebox(0,0)[lb]{$\bullet$}}
\put(1120,-1590){\makebox(0,0)[lb]{$\bullet$}}
\put(1120,-961){\makebox(0,0)[lb]{$\bullet$}}
\put(1120,-361){\makebox(0,0)[lb]{$\bullet$}}
\put(1120,239){\makebox(0,0)[lb]{$\bullet$}}
\put(1201,-1561){\vector( 1, 0){3600}}
\put(1201,-1561){\vector( 0, 1){2100}}
\put(2926,-436){\vector(-1,-1){450}}
\put(3526,-436){\vector(-1,-1){450}}
\put(2326,164){\vector(-1,-1){450}}
\put(2926,164){\vector(-1,-1){450}}
\put(3526,164){\vector(-1,-1){450}}
\put(4126,164){\vector(-1,-1){450}}
\put(4051,-436){\vector(-1,-1){450}}
\put(2926,764){\vector(-1,-1){450}}
\put(2326,764){\vector(-1,-1){450}}
\put(2326,-436){\vector(-1,-1){450}}
\put(4126,764){\vector(-1,-1){450}}
\put(3526,764){\vector(-1,-1){450}}
\put(1876,-436){\vector( 1,-1){450}}
\put(1876,164){\vector( 1,-1){450}}
\put(1876,764){\vector( 1,-1){450}}
\put(2476,764){\vector( 1,-1){450}}
\put(2476,164){\vector( 1,-1){450}}
\put(2476,-436){\vector( 1,-1){450}}
\put(3076,-436){\vector( 1,-1){450}}
\put(3076,164){\vector( 1,-1){450}}
\put(3076,764){\vector( 1,-1){450}}
\put(3676,764){\vector( 1,-1){450}}
\put(3676,164){\vector( 1,-1){450}}
\put(3676,-436){\vector( 1,-1){450}}
\put(4351,-1936){\makebox(0,0)[lb]{vertices}}
\put(526,389){\makebox(0,0)[lb]{edges}}
\put(1726,-436){\makebox(0,0)[lb]{$C^{1,2}$}}
\put(1726,164){\makebox(0,0)[lb]{$C^{1,3}$}}
\put(2326,-1036){\makebox(0,0)[lb]{$C^{2,1}$}}
\put(2326,-436){\makebox(0,0)[lb]{$C^{2,2}$}}
\put(2326,164){\makebox(0,0)[lb]{$C^{2,3}$}}
\put(2926,-1036){\makebox(0,0)[lb]{$C^{3,1}$}}
\put(2926,-436){\makebox(0,0)[lb]{$C^{3,2}$}}
\put(2926,164){\makebox(0,0)[lb]{$C^{3,3}$}}
\put(3526,-1036){\makebox(0,0)[lb]{$C^{4,1}$}}
\put(3526,-436){\makebox(0,0)[lb]{$C^{4,2}$}}
\put(3526,164){\makebox(0,0)[lb]{$C^{4,3}$}}
\put(1726,-1036){\makebox(0,0)[lb]{$C^{1,1}$}}
\end{picture}

On the picture South-West arrows $\swarrow$ represents $\dd$ differential, South-East arrows $\searrow$ represent $d$ differential. It is clear from the picture that the complex splits into a direct sum of two subcomplexes $C^{k,j}$ with $k+j$ even and $k+j$ odd. It is also clear that we are dealing with a total complex of a bicomplex. It implies that $C^{\cdot},\delta_{\alpha}$ are isomorphic for differnt values of $\alpha \neq 0,\infty$.

Introduce two diagonal complexes $C^{\cdot}_0$and $C^{\cdot}_1$:
\begin{align}
&C^{j}_0=\bigoplus_{k+j \equiv 0 \ mod \ 2}( C^{k,j})\notag \\
&C^{j}_1=\bigoplus_{k+j \equiv 1 \ mod \ 2}( C^{k,j})\notag
\end{align}

\begin{Prop}One has a following a selfevident description of cohomology of the complex \ref{E:oqx}
\begin{align}
&H^{even}=\bigoplus_{i \equiv 0 \ mod \ 2}H^i( C^{\cdot}_0) \oplus \bigoplus_{i \equiv 1 \ mod \ 2}H^i( C^{\cdot}_1)\notag \\
&H^{odd}=\bigoplus_{i \equiv 1 \ mod \ 2}H^i( C^{\cdot}_0) \oplus \bigoplus_{i \equiv 0 \ mod \ 2}H^i( C^{\cdot}_1)\notag \\
\end{align}
\end{Prop}

\subsection{Definition of the category ${\Bbb M}$.}

In this section we shall introduce a category  ${\Bbb M}$.. The category $\tilde {\Bbb M}$ is a very close kin of ${\Bbb M}$ will be introduced in the next section .  The  homology of  $\tilde {\Bbb M}$  is equal to  homology of the complex $C^{\cdot}$.

 The category  ${\Bbb M}$  splits into a disjoint  union of two subcategories ${\Bbb M}_{0,1}$.

 Objects of the category  ${\Bbb M}$ will be not nececerally connected {\bf fat graphs} of even euler characteristic (${\Bbb M}_{0}$-case) or odd Euler characterictic (${\Bbb M}_{1}$-case). As a reminder the adjective  ``fat''in our context  means that the germs of edges at  each vertex of a nondirected graph $\GG$ are cyclicly ordered.

We shall describe only generators for the set of morphisms. These are

a)Automorphisms of fat graphs.

b) Fix an edge $e \in \GG$ which is not a loop. Shrinking the edge $e$ to a point we are getting a new graph $\GG / e$. Then $e$ defines an element which we denote by the same letter $e \in Hom(\GG.\GG/e)$.

c) Another set of generator consists of pairs of germs of edges(we shall call them flags) at a vertex $v$ of a graph $\GG$.These flags must be not adjacent to eachother.  Such pair defines a morphism from graph $\GG$ to a graph $\GG'$ which we shalldescribe now.

 The graph $\GG'$ is equal to the graph $\GG$ outside a small neighborhood of the vertex $v$. The intersection of the graph $\GG$ with the  neighborhood is a graph $\GG_0$ with one vertex $v$ and a set   $F(v)$ of flags (we don't take into account the boundary point of flags). The morphism  excises  $\GG_0$ and replaces it by a new graph $\GG_0'$ with 3 connected components .

   The cyclic order on  $F(v)$  allowes us to split $F(v) \backslash  \{f_1,f_2\}$ into the union of two sets, which we denote  $F(v_1)$ and  $F(v_2)$. These sets carry a canonical cyclic order induced from $F(v)$. The sets  $F(v_1)$ and  $F(v_2)$ will be the flags at the vertices $v_1$, $v_2$ of the new graph $\GG_0'$. This describes two of the three connected components. The other component is the union of  $f_1,f_2$, By definition it has no interior vertices. 

It is clear that there is natural identification of the ends of $\GG_0$ and $\GG_0'$, so there is no question how to glue $\GG_0'$ into $\GG \backslash \GG_0$. 

 We illustrate our constriction of $\GG_0'$ from $\GG_0$ on the picture below. The little open circles are the vertices. The flags $f_1,f_2$ are marked by dots:

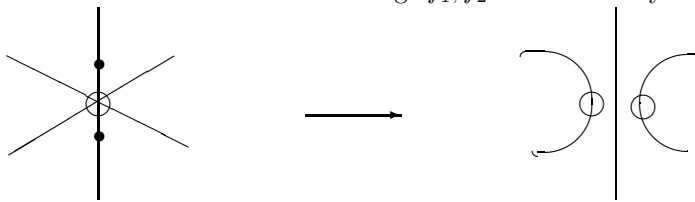
\setlength{\unitlength}{3947sp}%

\begin{picture}(4370,1224)(664,-973)
\thinlines
\put(4955,-61){\oval(142, 16)[tr]}
\put(4955,-357){\oval(608,608)[tl]}
\put(4951,-357){\oval(600,608)[bl]}
\put(4034,-661){\oval(116, 10)[bl]}
\put(4034,-349){\oval(634,634)[br]}
\put(4034,-349){\oval(634,634)[tr]}
\put(4034,-61){\oval(266, 58)[tl]}
\put(1215,-608){\makebox(0,0)[lb]{$\bullet $ }}
\put(1215,-158){\makebox(0,0)[lb]{$\bullet $ }}
\put(1250,-361){\circle{150}}
\put(4351,-361){\circle{150}}
\put(4673,-383){\circle{150}}
\put(1250,239){\line( 0,-1){1200}}
\put(676,-61){\line( 2,-1){1140}}
\put(1726,-61){\line(-5,-3){1036.765}}
\put(4501,239){\line( 0,-1){1200}}
\put(2551,-436){\vector( 1, 0){600}}
\end{picture}

There is an obvious definition of relations between morphisms which we skip.

\begin{Rem}
One can ``multiply'' objects in the category ${\Bbb M}$:
\begin{equation}\notag
{\Bbb M}_i \times {\Bbb M}_j \rightarrow {\Bbb M}_{i+j} \ mod \ 2
\end{equation}
 it is the union of graphs. 
\end{Rem}

\subsection{Definition of the categories $\tilde{\Bbb M}$}

We shall need to introduce some modification of the category ${\Bbb M}$ which we denote by  $\tilde{\Bbb M}$. As before It splits into a union of two connected components  $\tilde{\Bbb M}_{0,1}$.

The objects are the fat graphs with edges of two colors: black and white. 

The automorphisms of such graphs must preserve coloring.

The shrinking of edges are defined as usual.

In the third type of morphisms the marked flags must have the same color.

This category has a local system ${\cal O}$. We shall describe it in the next subsection, Now we can state the main 

\begin{Prop}

The cohomology of   $H^i(N(\tilde{\Bbb M}_{0,1}),{\cal O})$ are isomorphic to homology of the complex $C^i_{0,1}$ 
\end{Prop}

\subsection{Local system ${\cal O}$.}

By definition ${\cal O}(\GG)$ is a one dimensional vector space.It is a tensor product of certain number of even and odd elementary vector spaces.

Each black edge $e$ (it topologically is an open interval) provides an even group $H^1_{comp}(e,{\Bbb R})$.

At each vertex $v$ consider a set of all adjacent flags $F(v)$. The subset $F^w(v)$ is formed by white flags. For each $f \in F^w(v)$ define an odd vector space ${\Bbb R}_f$. Define a tensor product $S(v)=\Pi \bigotimes_{f \in F^w(v)} {\Bbb R}_f$. $\Pi$ stands for the parity change functor.

The  definition of the group ${\cal O}(\GG)$ is :
\begin{equation}\notag
 {\cal O}(\GG)= \bigotimes _{e \in black \ edges}H^1_{comp}(e,{\Bbb R})\otimes \bigotimes_{v \in vertices} S(v)
\end{equation}

\section{Questions}

One of the most intriguing questions is a relation of this construction to second quantized strings, Fukaya category on symmetric products $S^n(X)$, where $X$ is a Calabi-Yau manifold.

The other question is what is the relation of the category ${\Bbb M}$ to the moduli of {\bf all} curves? Such object was introduced in one of early attempt to define a nonperturbative string theory.

What is the relation to the mirror symmetry? There are no  quantum corrections on the holomorphic side. In what sence we can identify the extended strigy version of Fukaya category with apropriate modification of the category of vector bundles?  

The bialgebras are cocycles of deformations of commutative and \\ co-commutative polynomial Hopf algebras. Is it possible to find an actual deformation (quantization) and what is its physical interpretation(String theoretic,M theoretic,...)?


\begin{thebibliography}{9}
\bibitem{K}
M. Kontsevich, {\em Formal (non)-commutative symplectic geometry},
Progress in Mathematics, 107.  Birkhauser Boston, Inc.,
Boston, MA, 1993.
\bibitem{K1} M. Kontsevich, {\em  Homological algebra of mirror symmetry}, Proceedings of the International
Congress of Mathematicians, Vol. 1, 2 (Zurich, 1994), 120--139, Birkhauser, Basel, 1995.
\bibitem{M}M.Movshev {The definition of graph homology of algebras} To appear in Selecta Mathematica. New Series 
\bibitem{F}D.Freed {Five lecures on supersymmetry}  
\bibitem{Fuk}K. Fukaya, {Morse Homotopy, $A^\infty$-Category, and
Floer Homologies}, in  The Proceedings of the
1993 GARC Workshop on Geometry and Topology, H. J. Kim, ed.,
Seoul National University;
\end{thebibliography}
\end{document}